# Typability in partial groupoids

P.H. Rodenburg[1] and I. Bethke[2]

*Institute for Informatics, University of Amsterdam*

Adapting a claim of M. Kracht [K], we establish a characterization of the typable partial applicative algebras.

## Introduction

If you apply the copula ‘is’ to the adjective ‘dead’ you get the predicate ‘is dead’; and if you apply the noun ‘Socrates’ to this predicate, you get the sentence ‘Socrates is dead’. This is one way a grammarian could categorize words and groups of words and describe how they combine to form meaningful expressions. Now you could ask the opposite question: given certain combinations of words as *meaningful*, is there a grammatical categorization that explains these combinations? Kracht tried to answer this question (among other things) in [K], for a simple, abstract, grammatical model. We give a description of the details of this model below, state and refute a variant of Kracht’s theorem, and prove a more complex characterization along the same lines.

## Preliminaries

A *partial groupoid*, or more briefly, a *pargoid* [LE], consists of a set $A$ and a partial binary operation $\cdot$ (product) on $A$. If we think of the product operation as the *application* of a function to its argument, $\langle A, \cdot\rangle$ may be aptly referred to as a *partial applicative algebra*. As is the custom in combinatory logic, we tend to omit the operation symbol $\cdot$, and suppress parentheses assuming association to the left, writing $abc$ when we mean $(a \cdot b) \cdot c$. We write $ab\downarrow$ to express that $ab$ exists, and $ab\uparrow$ to express that it doesn’t. If we use $ab$ in a positive statement, such as that it belongs to some set (in particular, ‘$ab\uparrow$’ is a negative), we imply that it exists.

The polynomial operations of a pargoid $\mathbf{A}$ are the operations that can be constructed by composition from the product, projection and constant operations. The trivial polynomial operations are the ones that can be constructed without the product operation.

A congruence relation of a pargoid $\mathbf{A}$ is an equivalence relation $\rho$ of $A$ that respects the product in the sense that

$$a \equiv b \;\&\; c \equiv d \;\&\; ac\downarrow \;\&\; bd\downarrow \Rightarrow ac \equiv bd \;(\rho).$$

In particular, by $\varpi_{\mathbf{A}}$ we denote the congruence

[1] P.H.Rodenburg@uva.nl
[2] I.Bethke@uva.nl

$\{\langle a, b\rangle \mid$ for all unary polynomial operations $p$ of $\mathbf{A}$, $p(a)\downarrow \Leftrightarrow p(b)\downarrow\}$.

We observe that, relative to the definition of congruence that we just gave, $\varpi_{\mathbf{A}}$ is the Leibniz congruence $\Omega_{\mathbf{A}}(A)$ in the sense of Blok and Pigozzi [BP].

The quotient of a pargoid $\mathbf{A}$ over a congruence relation $\theta$ is $\mathbf{A}/\theta = \langle A/\theta, \cdot\rangle$.

A *type system* is an absolutely free algebra $\mathbf{T} = \langle T, \to\rangle$ with a single binary operation. The free generators of $\mathbf{T}$ are the *ground types*; the rest are *function types*. In our notation for function types we use association to the right: $\alpha \to \beta \to \gamma = \alpha \to (\beta \to \gamma)$. A subset $X$ of $T$ is *strict* if $\alpha \to \beta \in X$ only if $\alpha, \beta \in X$. A $\mathbf{T}$-*typed applicative algebra* is a pargoid $\mathbf{A} = \langle A, \cdot\rangle$ with an injective assignment $\alpha \mapsto A_\alpha$ of subsets of $A$ to the elements of a strict subset $S$ of $T$ such that, for $a, b \in A$,

1° $ab\downarrow \Leftrightarrow \exists \alpha, \beta \in S\, (a \in A_{\alpha\to\beta}\ \&\ b \in A_\alpha\ \&\ ab \in A_\beta)$,

2° $\{A_\alpha \mid \alpha \in S\}$ is a partition of $A$.

If such an assignment exists for $\mathbf{A}$, we say $\mathbf{A}$ is $\mathbf{T}$-*typable*. A pargoid is *typable* if for some type system $\mathbf{T}$ it is $\mathbf{T}$-typable. The type of an element of a $\mathbf{T}$-typed applicative algebra $\mathbf{A}$ is the unique $\alpha \in T$ such that $a \in A_\alpha$. The elements of $S$ are the *inhabited* types.

In a typable pargoid, no element applies to itself: by (1°), such an element $a$ should have a function type $\alpha \to \beta$, and also the antecedent type $\alpha$; so by (2°) and injectivity of the assignment $\xi \mapsto A_\xi$, $\alpha \to \beta = \alpha$, contradicting the absolute freedom of the type system.

**Lemma 1**. Let $\mathbf{A}$ be a typed applicative algebra. If $a, b \in A$ have the same type, then $a \equiv b$ $(\varpi_{\mathbf{A}})$.

**Proof**. By induction on unary polynomials, show that $p(a)$ and $p(b)$ have the same type if either one exists. ☒

### A claim and two counterexamples

Theorem 10 in [K] suggests the following characterization of typability for pargoids:

(*) A pargoid $\mathbf{A}$ is typable if and only if
(Tarski's Principle) $a \equiv c$ $(\varpi_{\mathbf{A}})$ if and only if there exists a nontrivial unary polynomial operation $p$ of $\mathbf{A}$ such that $p(a)\downarrow$ and $p(c)\downarrow$, and
(Well-Foundedness) for every $a$ there exist $n$ and $b_0,\ldots, b_n$ such that $ab_0\ldots b_n\uparrow$.

This is too simple; to see it is false, consider a pargoid with three elements $a$, $b$, $c$, and product specified by the table

| | $a$ | $b$ | $c$ |
|---|---|---|---|
| $a$ | $a$ | $\uparrow$ | $\uparrow$ |
| $b$ | $\uparrow$ | $\uparrow$ | $c$ |
| $c$ | $\uparrow$ | $\uparrow$ | $\uparrow$ |

Every element is in the domain of some nontrivial unary polynomial operation: $aa\downarrow$ and $bc\downarrow$, no nontrivial unary polynomial operation converges on more than one element, so $\varpi$ is the diagonal relation; and the product is well-founded since $ab$, $bb$ and $cb$ all diverge. However, $a$ cannot be given a type since it applies to itself.

This example indicates that the 'Well-Foundedness' condition is too weak. We shall formulate a better condition below.

The first condition, however, is problematic as well, on two counts. First, $a \equiv c$ $(\varpi)$ if all nontrivial unary polynomials diverge on $a$ and $c$. But, if the product operation of **A** is void, **A** is certainly typable. So the 'if and only if' should be *if*. And this will not be enough, for, second, consider the pargoid **A** specified by

| | $a$ | $b$ | $c$ | $ab$ | $cb$ | $d$ |
|---|---|---|---|---|---|---|
| $a$ | $\uparrow$ | $ab$ | $\uparrow$ | $\uparrow$ | $\uparrow$ | $\uparrow$ |
| $b$ | $\uparrow$ | $\uparrow$ | $\uparrow$ | $\uparrow$ | $\uparrow$ | $\uparrow$ |
| $c$ | $\uparrow$ | $cb$ | $\uparrow$ | $\uparrow$ | $\uparrow$ | $\uparrow$ |
| $ab$ | $\uparrow$ | $\uparrow$ | $\uparrow$ | $\uparrow$ | $\uparrow$ | $d$ |
| $cb$ | $\uparrow$ | $\uparrow$ | $\uparrow$ | $\uparrow$ | $\uparrow$ | $\uparrow$ |
| $d$ | $\uparrow$ | $\uparrow$ | $\uparrow$ | $\uparrow$ | $\uparrow$ | $\uparrow$ |

We have $a \not\equiv c$ $(\varpi_{\mathbf{A}})$, for $abd\downarrow$ whereas $cbd\uparrow$. But $ab$ and $cb$ both converge. So **A** should not be typable. But here is a typing of the elements:

$$a\colon \beta \to (\delta \to \delta);\quad b\colon \beta;\quad c\colon \beta \to \gamma;\quad ab\colon \delta \to \delta;\quad cb\colon \gamma;\quad d\colon \delta.$$

The problem is, that with 'nontrivial' we try to single out polynomial operations that really do something with their argument, but in $xb$, $x$ does something rather than that something is done with it.

**Definition**. The *definite* polynomial operations of **A** are the elements of the least class $P$ of nontrivial unary polynomial operations $p_1 \cdot p_2$ such that either $p_1 \in P$ or $p_2$ is nonconstant.

**Lemma 2**. If $p$ is a nonconstant unary polynomial operation that is indefinite, there are $b_1,\ldots, b_n \in A$ $(n \geq 0)$ such that $p(x) = xb_1\ldots b_n$.

**Proof**. By induction on polynomials. If $p$ is nonconstant and trivial, $p(x) = x$. If $p$ is indefinite and nontrivial, we must have $p(x) = p_1(x) \cdot p_2(x)$ with $p_2$ constant and $p_1$ nonconstant and indefinite. Use the induction hypothesis. ☒

**Characterization of typability**

Define, for elements $a$, $b$ of a pargoid $\mathbf{A} = \langle A, \cdot\rangle$,

$$b <_{\mathbf{A}} a \text{ if and only if: } ab{\downarrow} \text{ or } \exists c \in A.\ b = ac.$$

**Theorem**. A pargoid **A** is typable if and only if:

(i) for all $a$, $c \in A$, if there exists a definite polynomial operation $p$ of **A** such that $p(a){\downarrow}$ and $p(c){\downarrow}$, then $a \equiv c$ $(\varpi_{\mathbf{A}})$;

(ii) the relation $<_{\mathbf{A}}$ is well-founded.

**Proof**. ($\Rightarrow$) Assume a type system for **A**.
(i) If $p(a){\downarrow}$ and $p(c){\downarrow}$, for some definite polynomial operation $p$, then $a$ and $c$ have the same type. This is shown by induction on the construction of $p$. Since $p$ is definite, we may assume $p(x) = p_1(x) \cdot p_2(x)$, with $p_1$ definite or $p_2$ nonconstant. If $p_1$ is definite, by induction hypothesis $a$ and $c$ have the same type; likewise if $p_2$ is definite. If $p_2$ is indefinite and nonconstant, by Lemma 2, $p_2(x) = xb_1 \ldots b_n$, for certain $b_1, \ldots, b_n \in A$ $(n \geq 0)$ and all $x$. Let $\beta_1, \ldots, \beta_n$ be the respective types of $b_1, \ldots, b_n$. Any $x$ such that $p_2(x){\downarrow}$ must be of a type $\beta_1 \to \ldots \to \beta_n \to \xi$. If $p_1$ is constant, say $p_1(x) = d$, $d$ must be of some type $\xi \to \eta$. Since $d$ is fixed, this fixes $\xi$; hence $a$ and $c$ are of the same type. Otherwise we may assume that $p_1(x) = xd_1 \ldots d_m$, for certain $d_1, \ldots, d_m \in A$ $(n \geq 0)$ and all $x$. Let $\delta_1, \ldots, \delta_m$ be the respective types of $d_1, \ldots, d_m$. Now $x$ must also be of a type $\delta_1 \to \ldots \to \delta_m \to \xi \to \zeta$. By the freedom of the type system, $m < n$. So $\xi = \beta_i$, for some $i \leq n$; which again fixes the type of $x$. Now since the types of $a$ and $c$ are the same, by Lemma 1, $a \equiv c$ $(\varpi_{\mathbf{A}})$.
(ii) If $b <_{\mathbf{A}} a$, the type of $b$ is shorter than that of $a$.
($\Leftarrow$) Suppose **A** satisfies (i) and (ii). Let $\varpi$ be $\varpi_{\mathbf{A}}$, $<$ be $<_{\mathbf{A}}$. Define:

$$S_0 = \{\alpha \in A/\varpi \mid \exists a \in \alpha\ \forall b \in A\ ab{\uparrow}\}.$$

We let $S_0$ be the collection of ground types, and for $\alpha \in S_0$ put $A_\alpha = \alpha$. Function types are defined by

$$A_{\alpha\to\beta} = \{a \in A \mid \exists b \in A_\alpha\ ab \in A_\beta\}.$$

Let $a$ be a minimal element of $A - \mathrm{Y}_\alpha A_\alpha$. Then $a \notin \mathrm{Y}S_0$, so $ab{\downarrow}$ for some $b$. Since $b$, $ab < a$, there are $\alpha$, $\beta$ such that $b \in A_\alpha$ and $ab \in A_\beta$. Then $a \in A_{\alpha\to\beta}$. So every element has a type.

Now we prove by simultaneous induction on $<$:
$a \in A_\alpha \cap A_\beta$ implies $\alpha = \beta$;
$a \in A_\alpha$ implies $\forall b\ (b \equiv a\ (\varpi) \Rightarrow b \in A_\alpha)$.

Let $a$ be minimal among the elements that do not satisfy these conditions. Suppose $a \in \alpha \in S_0$. If $a \in A_\beta$, with $\alpha \neq \beta$, then $\beta \notin S_0$, for then $\alpha$ and $\beta$ would be distinct $\varpi$-congruence classes, and hence disjoint. So $\beta$ is a function type; say $ab{\downarrow}$. Let $c \in A_\alpha$ be such that $\forall d\ cd{\uparrow}$. Since $p(a){\downarrow}$, for $p(x) =$ $xb$, and $c \equiv a$ ($\varpi$), we have $cb = p(c)\downarrow$, quod non. The definition of the ground types ensures that $a$ satisfies the second condition.

Now suppose $a \in A_{\alpha\to\beta} \cap A_{\gamma\to\delta}$. Then there are $b \in A_\alpha$ and $c \in A_\gamma$ such that $ab \in A_\beta$ and $ac \in A_\delta$. So $p(b)\downarrow$ and $p(c)\downarrow$ for $p(x) = ax$. Hence by Tarski's Principle, $b \equiv c$ ($\varpi$). So by induction hypothesis, $\alpha = \gamma$. Since $\varpi$ is a congruence relation, a fortiori $ab \equiv ac$ ($\varpi$), so $\beta = \delta$. Hence $\alpha \to \beta = \gamma \to \delta$.

Finally, suppose $a \in A_{\alpha\to\beta}$ and $b \in a/\varpi$. Then for some $c \in A_\alpha$, $ac \in A_\beta$. Then for $p(x) = xc$: $p(a)\downarrow$, hence $p(b)\downarrow$, i.e. $bc\downarrow$. So $b \in A_{\alpha\to\gamma}$, with $bc \in A_\gamma$. Since $\varpi$ is a congruence relation, and $a \equiv b$ ($\varpi$), we have $ac \equiv bc$ ($\varpi$). Since $ac < a$, by induction hypothesis $\beta = \gamma$. So $\alpha \to \beta = \alpha \to \gamma$. ☒